# THE INFINITE IN SCIENCES AND ARTS

W. Mueckenheim

University of Applied Sciences Augsburg

*Actual infinity in its various forms is discussed, searched and not found.*

## INTRODUCTION

The infinite or, speaking precisely, the *actual infinite*[1], has been suspected in many domains. G. Cantor was convinced that the infinite exists in God, in many domains of nature and in mathematics. "Dementsprechend unterscheide ich ein 'Infinitum aeternum increatum sive Absolutum', das sich auf Gott und seine Attribute bezieht, und ein 'Infinitum creatum sive Transfinitum', das überall dort ausgesagt wird, wo in der Natura creata ein Aktual-Unendliches konstatiert werden muß, wie beispielsweise in Beziehung auf die, meiner festen Überzeugung nach, aktual-unendliche Zahl der geschaffenen Einzelwesen sowohl im Weltall wie auch schon auf unserer Erde und, aller Wahrscheinlichkeit nach, selbst in jedem noch so kleinen, ausgedehnten Teil des Raumes"[2] [1, p. 399]. In a letter to D. Hilbert he wrote about the applicaton of transfinite set theory "auf die *Naturwissenschaften*: Physik, Chemie, Mineralogie, Botanik, Zoologie, Anthropologie, Biologie, Physiologie, Medizin etc. ... Dazu kommen aber auch Anwendungen auf die sogenannten 'Geisteswissenschaften', die meines Erachtens als Naturwissenschaften aufzufassen sind, *denn auch der 'Geist' gehört mit zur Natur*"[3] [2, p. 459].

Cantor's first proofs of actual infinity had partially religious character - based on his firm belief in God and encouraged by the Holy Bible which says: "'Dominus regnabit in infinitum (aeternum) *et ultra*'. Ich meine dieses '*et ultra*' ist eine Andeutung dafür, dass es mit dem ω nicht sein Bewenden hat, sondern dass es auch *darüber hinaus* noch etwas giebt"[4] [2, p. 148]. "Ein Beweis geht vom Gottesbegriff aus und schließt zunächst aus der höchsten Vollkommenheit Gottes Wesens auf die Möglichkeit der Schöpfung eines Transfinitum ordinatum, sodann aus seiner Allgüte und Herrlichkeit auf die Notwendigkeit der tatsächlich erfolgten Schöpfung eines Transfinitum"[5] [1, p. 400].

---

[1] A potentially infinite quantity is always finite though not bounded, like a function. An actually infinite quantity is static, complete, i.e., *finished* and larger than any finite quantity of the same kind. Cantor calls it "finished infinite" (vollendet-unendlich) [1, p. 175], [2, p. 148, p. 181].

[2] Accordingly I distinguish an eternal uncreated infinity or absolutum which is due to God and his attributes, and a created infinity or transfinitum, which is has to be used wherever in the created nature an actual infinity has to be noticed, for example, with respect to, according to my firm conviction, the actually infinite number of created individuals, in the universe as well as on our earth and, most probably, even in every arbitrarily small extended piece of space.

[3] to the *natural sciences*: physics, chemistry, mineralogy, botany, zoology, anthropology, biology, physiology, medicine etc. ... In addition there are the applications to the so-called 'arts' which, in my view, also have to be considered natural sciences, *because also the 'mind' belongs to nature*.

In the following we will search high and low for evidence of actual infinity, be it as a reality or as an idea. With respect to *reality*, only physics including cosmology is concerned, while the question of actual infinity in the realm of theology (is it reality or idea?) is a difficult matter. Nevertheless, it belongs to the general topic and will be discussed. Actual infinity *as an idea* is to be suspected mainly in mathematics. In addition there is frequent mentioning of infinity in aesthetics and art but those effusive statements[6] rarely lay claim to cover the meaning of actual infinity and will not further be considered. Also subjective statements of historical scholars[7] as well as merely formal applications like the use of infinite vector spaces in theoretical physics do not belong to our topic.

## DOES THE INFINITELY SMALL EXIST IN REALITY?

Quarks are the smallest elementary particles presently known. Down to $10^{-19}$ m there is no structure detectable. Many physicists including the late W. Heisenberg are convinced that there is no deeper structure of matter. On the other hand, the experience with molecules, atoms, and elementary particles suggests that these physicists may be in error and that matter may be further divisible. However, it is not divisible in infinity. There is a clear-cut limit.

---

[4] 'The Lord rules in infinity (eternity) and beyond.' I think this 'and beyond' is a hint that ω is not the end of the tale.

[5] One proof is based on the notion of God. First, from the highest perfection of God, we infer the possibility of the creation of transfinity, then, from his all-grace and splendor, we infer the necessity that the creation of transfinity in fact has happened.

[6] "Silently one by one, in the infinite meadows of the heaven ..." (Henry W. Longfellow). "The Development of the aesthetics of the infinite" (Marjorie Hope Nicholson). "Above, infinite space, illimitable emptiness, with only the sun shining brazenly, eternally ..." (Sir Francis Chichester). " ... empty land stretching silently into infinity" (Alan Moorehead). "Infinity is where things happen that don't." (Anonymous). "I am painting the infinite" (Vincent van Gogh). "Towards the Infinite" (Joan Miró). [3]

[7] "Actual infinity does not exist in mathematics" (Aristoteles). "The number of points in a segment one ell long is its true measure." (R. Grosseteste). "Actual infinity exists in number, time and quantity" (J. Baconthorpe). "Je suis tellement pour l'infini actuel" (G. Leibniz). [4]

Lengths which are too small to be handled by material meter sticks can be measured in terms of wavelengths $\lambda$ of electromagnetic waves, for instance.

$$\lambda = c/\nu \qquad (c = 3*10^8 \text{ m/s})$$

The frequency $\nu$ is given by the energy $E$ of the photon

$$\nu = E/h \qquad (h = 6{,}6*10^{-34} \text{ Js})$$

and a photon cannot contain more than all the energy of the universe

$$E = mc^2$$

which has a mass of about $m = 5*10^{55}$ g (including the dark matter). This yields the complete energy $E = 5*10^{69}$ J. So the unsurpassable minimal length is $4*10^{-95}$ m.

## DOES THE INFINITELY LARGE EXIST IN REALITY?

Modern cosmology teaches us that the universe has a beginning and is finite. But even if we do not trust in this wisdom, we know that theory of relativity is as correct as human knowledge can be. According to relativity theory, the accessible part of the universe is a sphere of $50*10^9$ light years radius containing a volume of $10^{80}$ m$^3$. (This sphere is growing with time but will remain finite forever.) "Warp" propulsion, "worm hole" traffic, and other science fiction (and scientific fiction) does not work without time reversal. Therefore it will remain impossible to leave (and to know more than) this finite sphere. Modern quantum mechanics has taught us that entities which are non-measurable in principle, do not exist. Therefore, also an upper bound (which is certainly not the supremum) of $10^{365}$ for the number of elementary spatial cells in the universe can be calculated from the minimal length estimated above.

## HOW LONG LASTS ETERNITY?

We know that the universe has been expanding for about $14*10^9$ years. This process will probably continue in eternity. According to newer astronomic results the over-all energy density is very close to zero, suggesting an Euclidean space. Eternity, however, will never be completed. So time like space has a potentially infinite character. Both are of unbounded size but always finite.

Will intelligent creatures survive in eternity? One constraint is the limited supply of free energy which is necessary for any life. This problem could be solved, however, by living for shorter and shorter intervals with long hibernation phases. By means of series like $1/2 + 1/4 + 1/8 + ... = 1$, a limited amount of energy could then last for ever - and with it intelligent life [5].

Alas, there is the risk of sudden death of a creature by an accident. If we assume that in our civilisation one out of 200 lives ends by an accident, then we can calculate that the risk to die by an accident during this very minute is roughly $10^{-10}$. The risk that 6,000,000,000 people will die during this minute is then $10^{-60,000,000,000}$ (not taking into account cosmic catastrophes, epidemic diseases etc.). This is a very small but positive probability. And even if in future the risk of accidental death can be significantly reduced while the population may be enormously increased[8], the risk of a sudden end of all life within a short time interval is not zero and, therefore, will occur before eternity is finished [6, p. 135].

So, after a while, nobody will be present to measure time - and it may well be asked if an entity does exist which in principle cannot be measured.

## IS THE HEAVEN INFINITE?

The existence of a creator has become more and more improbable in history of mankind. Copernicus, Bruno, and Darwin contributed to remove mankind from the centre of the universe, the position chosen by God for his creatures. The development of the character of God himself closely reflects the social development of human societies. G. C. Lichtenberg observed: "God created man according to his image? That means probably, man created God according to his."

The discovery of foreign cultures, in America and Australia, showed that people there had not been informed in advance about the God of Jews, Christians or Moslems - a highly unfair state of affairs in case belief in this God was advantageous before or after death.

The results of neurology and cerebral surgery show that characteristic traits and behaviour usually retraced to the human soul can be arbitrarily manipulated by electric currents, drugs, or surgery while an immortal soul cannot be localized.

Of course it is impossible to prove or to disprove the existence of one or more Gods, but it is easy to disprove the absolutum, as Cantor called it, i.e., the infinity of every property of a God. Medieval scholastics already asked, whether God could make a stone that heavy that he himself was incapable of lifting it. God cannot know the complete future unless the universe is deterministic. But in this case, there could be no free will and no living creature could prove itself suitable or unsuitable to enter paradise or hell - and the whole creation was meaningless.

Therefore, actual infinity, as being inherent to theological items, cannot be excluded but is at best problematic.

## ARE THERE INFINITE SETS IN MATHEMATICS?

Set theory contains an axiom postulating the existence of an actually infinite set, and its most important theorem states that there are even different infinities. In particular the cardinal number of the infinite set of real numbers is larger than the cardinal number of the "smallest" infinite set, namely the set of natural numbers. This result has raised differing opinions. While most mathematicians were pleased about Cantor's theory "diese erscheint mir als die bewundernswerteste Blüte mathemati-schen Geistes und überhaupt eine der höchsten Leistungen rein verstandesmäßiger menschlicher Tätigkeit"[9] [7, p. 167], "aus dem Paradies, das Cantor für uns geschaffen, soll uns niemand vertreiben können"[10] [7, p. 170], others remained sceptic. "In der Renaissance, besonders bei Bruno, überträgt sich die aktuale Unendlichkeit von Gott auf die Welt. Die endlichen Weltmodelle der gegenwärtigen Naturwissenschaft zeigen deutlich, wie diese Herrschaft des Gedankens einer aktualen Unendlichkeit mit der klassischen (neuzeitlichen) Physik zu Ende gegangen ist. Befremdlich wirkt dem gegenüber die Einbeziehung des Aktual-Unendlichen in die Mathematik, die explizit erst gegen Ende des vorigen Jahrhunderts mit G. Cantor begann. Im geistigen Gesamtbilde unseres Jahrhunderts ... wirkt das aktual Unendliche geradezu anachronistisch"[11] [8]. "Infinite totalities do not exist in any sense of the word (i.e., either really or ideally). More precisely, any mention, or purported mention, of infinite totalities is, literally, meaningless" [9].

---

[8] but at most to $10^{78}$ individuals because this is the number of hydrogen atoms in the universe, and life without at least one atom seems impossible.

In fact, there is some evidence that actual infinity, at least if understood as a number larger than every natural number, leads to problems. One simple observation is that every initial segment of even positive integers has a cardinal number

$$|\{2, 4, 6, ..., 2n\}| = n$$

which is less than some of the elements. There is no reason to expect that this discrepancy will disappear for the union of all initial segments, because it increases with increasing size of the segment. But this union is the set of all even positive integers (because there is no positive integer outside of the union). Its cardinal number should be larger than any element of the set.

The most important theorem of set theory states that the set of real numbers of the interval [0, 1] has a greater cardinal number than the set of natural numbers, i.e., the continuum [0, 1] is uncountable.

The infinite binary tree $T$, consisting of nodes, namely the numerals 0 and 1, contains the binary representations of all real numbers of the interval [0, 1] in form of infinite paths, i.e., sequences of nodes, starting at the root node and running through the levels of the tree, guided by the edges, as indicated below. Some real numbers, like 0.1000... = 0.0111...., are present even in duplicate.

---

[9] This is, in my opinion, the most admirable blossom of the mathematical mind and altogether one of the outstanding achievements of purely intellectual human activity.

[10] No one shall expel us from the paradise Cantor created for us.

[11] During the renaissance, particularly with Bruno, actual infinity transfers from God to the world. The finite world models of contemporary science clearly show how this power of the idea of actual infinity has ceased with classical (modern) physics. Under this aspect, the inclusion of actual infinity into mathematics, which explicitly started with G. Cantor only towards the end of the last century, seems displeasing. Within the intellectual overall picture of our century ... actual infinity brings about an impression of anachronism.

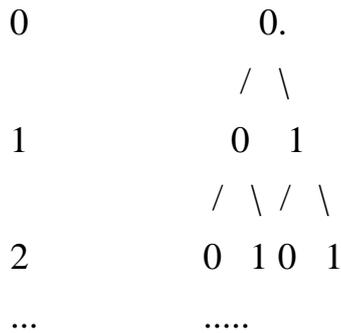

```
     0              0.
                   / \
     1            0   1
                 / \ / \
     2          0  1 0  1
    ...           .....
```

Every *finite* binary tree $T_n$ contains less paths than nodes. Down to level *n* there are $2^{n+1} - 1$ nodes but only $2^n$ path. In the infinite binary tree *T* the set of nodes remains countable. But the set of paths must be uncountable, because for each real number of the interval [0, 1] there is at least one path in *T* representing it.

The union $T_\infty$ of *all finite* binary trees *covers all levels* enumerated by natural numbers. With respect to nodes and edges it is identical with the infinite binary tree

$$T_\infty = T.$$

According to set theory (including the axiom of choice) a countable union of countable sets is a countable set. The set of paths in the union tree $T_\infty$ is merely a countable union of *finite* sets, and, therefore, $T_\infty$ contains only a countable set of paths. But does $T_\infty$ contain only finite paths?

An *index* denotes the level to which a node belongs. The union of all indexes of nodes of finite paths is the union of all initial segments of natural numbers

$$\{1\} \cup \{1, 2\} \cup \{1, 2, 3\} \cup ... \cup \{1, 2, 3, ..., n\}... \cup ... = \{1, 2, 3, ...\}.$$

This is also the set of all last elements of the finite segments, i.e., it is the set of all natural numbers. This is the set of *all* indexes - there is no one left out. With respect to this observation we examine, for instance, all finite paths of the tree $T_\infty$ which always turn right: 0.1, 0.11, 0.111, ... If considered as sets of nodes, their union is the *infinite* path representing the real number 0.111... = 1. From this we can conclude that also *every* other infinite path belongs to the union $T_\infty$ of all finite trees.

The trees $T_\infty$ and *T* are identical with respect to all nodes, all edges, and all paths (which would already have been implied by the identity of nodes and edges). But the set of all paths is countable in the tree $T_\infty$ and uncountable in the same tree *T*.

Such problems of set theory can be avoided if we replace the fundamental concept of a one-to-one correspondence or bijection (leading to different transfinite cardinal numbers) by the concept of *intercession* [6, p. 116] briefly outline below.

Two infinite sets, *A* and *B*, *intercede* (each other) if they *can be* put in an intercession, i.e., if they can be ordered such that between two elements of *A* there is at least one element of *B* and, vice versa, between two elements of *B* there is at least one element of *A*.

The intercession includes Cantor's definition of equivalent (or equipotent) sets. Two equivalent sets always intercede each other, i.e., they can always be put in an intercession. (The intercession of sets with nonempty intersection, e.g., the intercession of a set with itself, requires the distinction of identical elements.[12]) The intercession is an equivalence relation, alas it is not as exciting as the bijection. *All infinite sets (like the integers, the rationals, and the reals) belong, under this relation, to one and the same equivalence class*. The sets of rational numbers and irrational numbers, for instance, intercede already in their natural order. There is no playing ground for building hierarchies upon hierarchies of infinities, for accessing inaccessible numbers, and for finishing the infinite. Every set which is not finite, is simply infinite, namely potentially infinite.

This point of view is in much better accordance with reality than the belief in actual infinity which has no application outside of mathematics at all [10, 11]. There are no numbers to be taken from an infinite set or from a Platonist shelf. "Die Zahlen sind freie Schöpfungen des menschlichen Geistes"[13] [12, p. III]. "A construction does not exist until it is made; when something new is made, it is something new and not a selection from a pre-existing collection" [13, p. 2]. When the objects of discussion are linguistic entities ... then that collection of entities may vary as a result of discussion about them. A consequence of this is that the 'natural numbers' of today are not the same as the 'natural numbers' of yesterday" [14, p. 478]. And we may add, what differs is not constant, i.e, it is not a constant set.

## CONCLUSION

We have investigated actual infinity in all domains which reasonably could be suspected to sustain it. But the result is such that we can give the last word to D. Hilbert, who, at the end of his famous paper praising actual infinity of set theory, comes to a plain and astonishing conclusion:

"Zuletzt wollen wir wieder unseres eigentlichen Themas gedenken und über das Unendliche das Fazit aus allen unseren Überlegungen ziehen: Das Gesamtergebnis ist dann: das Unendliche findet sich nirgends realisiert; es ist weder in der Natur vorhanden, noch als Grundlage in unserem verstandesmäßigen Denken zulässig - eine bemerkenswerte Harmonie zwischen Sein und Denken"[14] [7, 190].

---

[12] As an example an intercession of the set of positive integers and the set of even positive integers is given by {1, 2, 3, ...} ∪ {2', 4', 6', ...} = {1, 2', 2, 4', 3, 6', ...}.

[13] Numbers are free creations of the human mind.

[14] Finally, let us return to our original topic, and let us draw the conclusion from all our reflections on the infinite. The overall result is then: The infinite is nowhere realized. Neither is it present in nature nor is it admissible as a foundation of our rational thinking - a remarkable harmony between being and thinking.